\documentclass[11pt]{amsart}
\usepackage{amsmath,amsthm,amssymb}
\usepackage[arrow,matrix]{xy}
\usepackage{hyperref,a4wide}
\usepackage{pstricks,pst-node,multido}

\newtheorem{thm}{Theorem}[section]

\newtheorem{lem}{Lemma}[section]
\makeatletter
\renewcommand*{\c@lem}{\c@thm}
\makeatother

\newtheorem{prop}{Proposition}[section]
\makeatletter
\renewcommand*{\c@prop}{\c@thm}
\makeatother

\newtheorem{cor}{Corollary}[section]
\makeatletter
\renewcommand*{\c@cor}{\c@thm}
\makeatother

\theoremstyle{definition}

\newtheorem{example}{Example}[section]
\makeatletter
\renewcommand*{\c@example}{\c@thm}
\makeatother

\newtheorem{defn}{Definition}[section]
\makeatletter
\renewcommand*{\c@defn}{\c@thm}
\makeatother

\makeatletter
\renewcommand*{\c@conj}{\c@thm}
\makeatother

\makeatletter
\renewcommand*{\c@question}{\c@thm}
\makeatother

\makeatletter
\renewcommand*{\c@problem}{\c@thm}
\makeatother

\theoremstyle{remark}
\newtheorem{remark}{Remark}[section]
\makeatletter
\renewcommand*{\c@remark}{\c@thm}
\makeatother

\DeclareMathOperator{\nc}{cc}
\def\k{\mathbf{k}}

\begin{document}
\title[A combinatorial proof on Betti numbers] {A combinatorial proof
  of a formula for Betti numbers of a stacked polytope}

\author{Suyoung Choi}
\address{Department of Mathematical Sciences, KAIST, 335 Gwahangno, Yuseong-gu, Daejeon 305-701, Republic of Korea}
\email{choisy@kaist.ac.kr}

\author{Jang Soo Kim}
\address{Department of Mathematical Sciences, KAIST, 335 Gwahangno, Yuseong-gu, Daejeon 305-701, Republic of Korea}
\email{jskim@kaist.ac.kr}

\date{\today}
\maketitle

\psset{dimen=inner,linewidth=0.5pt}

\begin{abstract}
For a simplicial complex $\Delta$, the graded Betti number
$\beta_{i,j}(\k[\Delta])$ of the Stanley-Reisner ring $\k[\Delta]$
over a field $\k$ has a combinatorial interpretation due to Hochster.
Terai and Hibi showed that if $\Delta$ is the boundary complex of a
$d$-dimensional stacked polytope with $n$ vertices for $d\geq3$, then
$\beta_{k-1,k}(\k[\Delta])=(k-1)\binom{n-d}{k}$. We prove this
combinatorially.
\end{abstract}


\section{Introduction}

A \emph{simplicial complex} $\Delta$ on a finite set $V$ is a
collection of subsets of $V$ satisfying
\begin{enumerate}
\item if $v\in V$ then $\{v\}\in \Delta$,
\item if $F\in \Delta$ and $F'\subset F$, then $F'\in \Delta$.
\end{enumerate}
Each element $F\in\Delta$ is called a \emph{face} of $\Delta$.  The
\emph{dimension} of $F$ is defined by $\dim(F)=|F|-1$.  The
\emph{dimension} of $\Delta$ is defined by
$\dim(\Delta)=\max\{\dim(F):F\in\Delta\}$.  For a subset $W\subset V$,
let $\Delta_W$ denote the simplicial complex $\{F\cap W:F\in\Delta\}$
on $W$.

Let $\Delta$ be a simplicial complex on $V$.  Two elements $v,u\in V$
are said to be \emph{connected} if there is a sequence of vertices
$v=u_0, u_1,\ldots, u_r=u$ such that $\{u_i,u_{i+1}\}\in \Delta$ for
all $i=0,1,\ldots,r-1$.  A \emph{connected component} $C$ of $\Delta$
is a maximal nonempty subset of $V$ such that every two elements of
$C$ are connected.  

Let $V=\{x_1,x_2,\ldots,x_n\}$ and let $R$ be the polynomial ring
$\k[x_1,\ldots,x_n]$ over a fixed field $\k$.  Then $R$ is a graded
ring with the standard grading $R=\oplus_{i\geq0}R_i$.  Let
$R(-j)=\oplus_{i\geq0}(R(-j))_i$ be the graded module over $R$ with
$(R(-j))_i=R_{j+i}$.  The \emph{Stanley-Reisner ring} $\k[\Delta]$ of
$\Delta$ over $\k$ is defined to be $R/I_{\Delta}$, where $I_{\Delta}$
is the ideal of $R$ generated by the monomials $x_{i_1}x_{i_2}\cdots
x_{i_r}$ such that $\{x_{i_1},x_{i_2},\ldots,x_{i_r}\}\not\in\Delta$.
A \emph{finite free resolution} of $\k[\Delta]$ is an exact sequence
\begin{equation}\label{eq:ffr}
\xymatrix{ 0 \ar[r]& F_r \ar[r]^{\phi_r}& F_{r-1} \ar[r]^{\phi_{r-1}}&
  \cdots \ar[r]^{\phi_{2}} & F_1\ar[r]^{\phi_1}& F_0 \ar[r]^{\phi_0}&
  \k[\Delta] \ar[r] & 0},
\end{equation}
where $F_i=\oplus_{j\geq0} R(-j)^{\beta_{i,j}}$ and each $\phi_i$ is
degree-preserving. A finite free resolution \eqref{eq:ffr} is
\emph{minimal} if each $\beta_{i,j}$ is smallest possible.  There is a
minimal finite free resolution of $\k[\Delta]$ and it is unique up to
isomorphism. If \eqref{eq:ffr} is minimal, then the \emph{$(i,j)$-th
  graded Betti number} $\beta_{i,j}(\k[\Delta])$ of $\k[\Delta]$ is
defined to be $\beta_{i,j}(\k[\Delta])=\beta_{i,j}$.  Hochster's
theorem says
$$\beta_{i,j}(\k[\Delta])=\sum_{\substack{W\subset V\\ |W|=j}}
\dim_{\k}\widetilde H_{j-i-1}(\Delta_W;\k).$$ We refer the reader to
\cite{Bruns1993,Stanley1996} for the details of Betti numbers and
Hochster's theorem.  Since $\dim_{\k}\widetilde H_{0}(\Delta_W;\k)$ is
the number of connected components of $\Delta_W$ minus $1$, we can
interpret $\beta_{i-1,i}(\k[\Delta])$ in a purely combinatorial way.

\begin{defn}
Let $\Delta$ be a simplicial complex on a finite nonempty set $V$. Let
$k$ be a nonnegative integer. The \emph{$k$-th special graded Betti
  number} $b_k(\Delta)$ of $\Delta$ is defined to be
\begin{equation}\label{eq:bk}
b_k(\Delta)=\sum_{\substack{W\subset V\\|W|=k}}
\left(\nc(\Delta_W)-1\right),
\end{equation}
where $\nc(\Delta_W)$ denotes the number of connected components of
$\Delta_W$.
\end{defn}
Note that since there is no connected component in $\Delta_{\emptyset}
=\{\emptyset\}$, we have $b_0(\Delta)=-1$.  If $k>|V|$, then
$b_k(\Delta)=0$ because there is nothing in the sum in
\eqref{eq:bk}. Thus we have
$$b_k(\Delta)=\left\{
\begin{array}{ll}
\beta_{k-1,k}(\k[\Delta]), & \mbox{if $k\geq1$,}\\
-1, & \mbox{if $k=0$}.\\
\end{array}
\right.$$

We refer the reader to \cite{Ziegler1995} for the basic notions of
convex polytopes. Let $P$ be a simplicial polytope with vertex set
$V$. The \emph{boundary complex} $\Delta(P)$ is the simplicial complex
$\Delta$ on $V$ such that $F\in\Delta$ for some $F\subset V$ if and
only if $F\ne V$ and the convex hull of $F$ is a face of $P$.  Note
that if the dimension of $P$ is $d$, then $\dim(\Delta(P))=d-1$.

For a $d$-dimensional simplicial polytope $P$, we can attach a
$d$-dimensional simplex to a facet of $P$.  A \emph{stacked polytope}
is a simplicial polytope obtained in this way starting with a
$d$-dimensional simplex.

Let $P$ be a $d$-dimensional stacked polytope with $n$ vertices.  Hibi
and Terai \cite{Terai1997} showed that $\beta_{i,j}(\k[\Delta(P)])=0$
unless $i=j-1$ or $i=j-d+1$. Since
$\beta_{i-1,i}(\k[\Delta(P)])=\beta_{n-i-d+1,n-i}(\k[\Delta(P)])$, it
is sufficient to determine $\beta_{i-1,i}(\k[\Delta(P)])$ to find all
$\beta_{i,j}(\k[\Delta(P)])$.  In the same paper, they found the
following formula for $\beta_{k-1,k}(\k[\Delta(P)])$:
\begin{equation}\label{eq:hibi}
\beta_{k-1,k}(\k[\Delta(P)])=(k-1)\binom{n-d}{k}.
\end{equation}
Herzog and Li Marzi \cite{Herzog1999} gave another proof of
\eqref{eq:hibi}.

The main purpose of this paper is to prove \eqref{eq:hibi}
combinatorially.  In the meanwhile, we get as corollaries the results
of Bruns and Hibi \cite{Bruns1998} : a formula of $b_k(\Delta)$ if
$\Delta$ is a tree (or a cycle) considered as a $1$-dimensional
simplicial complex.

\section{Definition of $t$-connected sum}
In this section we define a $t$-connected sum of simplicial complexes,
which gives another equivalent definition of the boundary complex of a
stacked polytope. See \cite{Buchstaber2002} for the details of
connected sums. And then, we extend the definition of $t$-connected
sum to graphs, which has less restrictions on the construction.  Every
graph in this paper is simple.

\subsection{A $t$-connected sum of simplicial complexes}
Let $V$ and $V'$ be finite sets. A \emph{relabeling} is a bijection
$\sigma:V\rightarrow V'$.  If $\Delta$ is a simplicial complex on $V$,
then $\sigma(\Delta) =\{\sigma(F):F\in\Delta\}$ is a simplicial
complex on $V'$.

\begin{defn}
Let $\Delta_1$ and $\Delta_2$ be simplicial complexes on $V_1$ and
$V_2$ respectively.  Let $F_1\in\Delta_1$ and $F_2\in\Delta_2$ be
maximal faces with $|F_1|=|F_2|$. Let $V_2'$ be a finite set and
$\sigma:V_2\rightarrow V_2'$ a relabeling such that $V_1\cap V_2'=F_1$
and $\sigma(F_2)=F_1$. Then the \emph{connected sum}
$\Delta_1\#^{F_1,F_2}_{\sigma}\Delta_2$ of $\Delta_1$ and $\Delta_2$
with respect to $(F_1,F_2,\sigma)$ is the simplicial complex
$(\Delta_1\cup\sigma(\Delta_2))\setminus\{F_1\}$ on $V_1\cup V_2'$.
If $\Delta=\Delta_1\#^{F_1,F_2}_{\sigma}\Delta_2$ and $|F_1|=|F_2|=t$,
then we say that $\Delta$ is a \emph{$t$-connected sum} of $\Delta_1$
and $\Delta_2$.
\end{defn}
Note that if $\Delta_1$ and $\Delta_2$ are $(d-1)$-dimensional pure
simplicial complexes, i.e. the dimension of each maximal face is
$d-1$, then we can only define a $d$-connected sum of them.

Let $\Delta_1,\Delta_2,\ldots,\Delta_n$ be simplicial complexes.  A
simplicial complex $\Delta$ is said to be a $t$-connected sum of
$\Delta_1,\Delta_2,\ldots,\Delta_n$ if there is a sequence of
simplicial complexes $\Delta_1',\Delta_2',\ldots,\Delta_n'$ such that
$\Delta_1'=\Delta_1$, $\Delta_i'$ is a $t$-connected sum of
$\Delta_{i-1}'$ and $\Delta_i$ for $i=2,3,\ldots,n$, and
$\Delta_n'=\Delta$.

\subsection{A $t$-connected sum of graphs}

Let $G$ be a graph with vertex set $V$ and edge set $E$. Let $W\subset
V$.  Then the \emph{induced subgraph} $G|_W$ of $G$ with respect to
$W$ is the graph with vertex set $W$ and edge set $\{\{x,y\}\in E: x,
y\in W\}$.  Let $$b_k(G)=\sum_{\substack{W\subset V\\|W|=k}}
\left(\nc(G|_W)-1\right),$$ where $\nc(G|_W)$ denotes the number of
connected components of $G|_W$. 

Let $\Delta$ be a simplicial complex on $V$.  The \emph{$1$-skeleton}
$G(\Delta)$ of $\Delta$ is the graph with vertex set $V$ and edge set
$E=\{F\in \Delta: |F|=2\}$. By definition, the connected components of
$\Delta_W$ and $G(\Delta)|_W$ are identical for all $W\subset V$. Thus
$b_k(\Delta)=b_k(G(\Delta))$.

Now we define a $t$-connected sum of two graphs.

\begin{defn}
Let $G_1$ and $G_2$ be graphs with vertex sets $V_1$ and $V_2$, and
edge sets $E_1$ and $E_2$ respectively.  Let $F_1\subset V_1$ and
$F_2\subset V_2$ be sets of vertices such that $|F_1|=|F_2|$, and
$G_1|_{F_1}$ and $G_2|_{F_2}$ are complete graphs.  Let $V_2'$ be a
finite set and $\sigma:V_2\rightarrow V_2'$ a relabeling such that
$V_1\cap V_2'=F_1$ and $\sigma(F_2)=F_1$. Then the \emph{connected
  sum} $G_1\#^{F_1,F_2}_{\sigma}G_2$ of $G_1$ and $G_2$ with respect
to $(F_1,F_2,\sigma)$ is the graph with vertex set $V_1\cup V_2'$ and
edge set $E_1\cup\sigma(E_2)$, where
$\sigma(E_2)=\{\{\sigma(x),\sigma(y)\}:\{x,y\}\in E_2\}$.  If
$G=G_1\#^{F_1,F_2}_{\sigma}G_2$ and $|F_1|=|F_2|=t$, then we say that
$G$ is a \emph{$t$-connected sum} of $G_1$ and $G_2$.
\end{defn}
Note that in contrary to the definition of $t$-connected sum of
simplicial complexes, it is not required that $F_1$ and $F_2$ are
maximal, and we do not remove any element in $E_1\cup\sigma(E_2)$.  We
define a $t$-connected sum of $G_1,G_2,\ldots,G_n$ as we did for
simplicial complexes.

It is easy to see that, if $|F_1|=|F_2|\geq3$ then
$G(\Delta_1\#^{F_1,F_2}_{\sigma}\Delta_2)
=G(\Delta_1)\#^{F_1,F_2}_{\sigma}G(\Delta_2)$. 
Thus we get the following proposition.

\begin{prop}\label{thm:tcon}
For $t\geq3$, if $\Delta$ is a $t$-connected sum of
$\Delta_1,\Delta_2,\ldots,\Delta_n$, then $G(\Delta)$ is a
$t$-connected sum of $G(\Delta_1),G(\Delta_2),\ldots,G(\Delta_n)$.
\end{prop}

Note that \autoref{thm:tcon} is not true if $t=2$ as the following
example shows.

\def\vput[#1](#2)#3{\pnode(#2){#3} \uput{.1}[#1](#2){$#3$}}
\def\edge#1#2{\ncline{#1}{#2}}

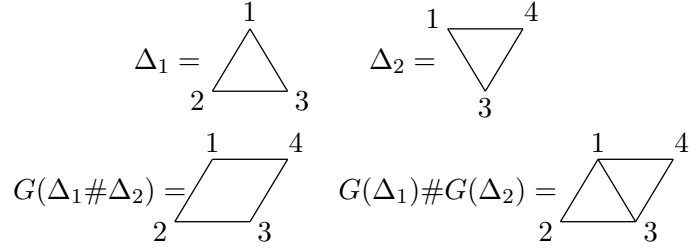
\begin{figure}
\begin{pspicture}(-1,0)(1,1)
\rput(-.6,.4){$\Delta_1=$}
\vput[210](0,0)2 \vput[90](.5,.83)1 \vput[330](1,0)3
\edge12 \edge13 \edge23
\end{pspicture}
\begin{pspicture}(-2,0)(1,1)
\rput(-.6,.4){$\Delta_2=$}
\vput[150](0,.83)1 \vput[270](.5,0)3 \vput[60](1,.83)4
\edge14 \edge13 \edge43
\end{pspicture}

\begin{pspicture}(-2,0)(1,1.7)
\rput(-1,.4){$G(\Delta_1\#\Delta_2)=$}
\vput[210](0,0)2 \vput[90](.5,.83)1 \vput[330](1,0)3 \vput[60](1.5,0.83)4
\edge12 \edge23 \edge43 \edge14
\end{pspicture}
\begin{pspicture}(-4,0)(1,1.7)
\rput(-1.5,.4){$G(\Delta_1)\# G(\Delta_2)=$}
\vput[210](0,0)2 \vput[90](.5,.83)1 \vput[330](1,0)3 \vput[60](1.5,0.83)4
\edge12 \edge23 \edge43 \edge14 \edge13
\end{pspicture}
\caption{The $1$-skeleton of a $2$-connected sum of $\Delta_1$ and
  $\Delta_2$ is not a $2$-connected sum of $G(\Delta_1)$ and
  $G(\Delta_2)$.}\label{fig:counter}
\end{figure}

\begin{example}
Let $\Delta_1=\{12,23,13\}$ and $\Delta_2=\{13,34,14\}$ be simplicial
complexes on $V_1=\{1,2,3\}$ and $V_2=\{1,3,4\}$. Here $12$ means the
set $\{1,2\}$. Let $F_1=F_2=\{1,3\}$ and let $\sigma$ be the identity
map from $V_2$ to itself. Then the edge set of
$G(\Delta_1\#^{F_1,F_2}_{\sigma}\Delta_2)$ is $\{12,23,34,14\}$, but
the edge set of $G(\Delta_1)\#^{F_1,F_2}_{\sigma}G(\Delta_2)$ is
$\{12,23,34,14,13\}$. See \autoref{fig:counter}.
\end{example}

\section{Main results}

In this section we find a formula of $b_k(G)$ for a graph $G$ which is
a $t$-connected sum of two graphs. To do this let us introduce the
following notation. For a graph $G$ with vertex set $V$, let $$c_k(G)
=\sum_{\substack{W\subset V\\|W|=k}}\nc(G|_W).$$ Note that
$c_k(G)=b_k(G)+\binom{|V|}{k}$.

\begin{lem}\label{lem:connectedsum}
Let $G_1$ and $G_2$ be graphs with $n_1$ and $n_2$ vertices
respectively. Let $t$ be a positive integer and let $G$ be a
$t$-connected sum of $G_1$ and $G_2$. Then
\begin{align*}
c_k (G) =& \sum_{i=0}^k \left( c_i (G_1) \binom{n_2 - t}{k-i} + c_i
(G_2) \binom{n_1 - t}{k-i}\right) \\
& - \binom{n_1 + n_2 -t}{k} + \binom{n_1 + n_2 -2t}{k}.  
\end{align*}
\end{lem}
\begin{proof}
Let $V_1$ (resp. $V_2$) be the vertex set of $G_1$ (resp. $G_2$).  We
have $G=G_1\#^{F_1,F_2}_{\sigma}G_2$ for some $F_1\subset V_1$,
$F_2\subset V_2$, a vertex set $V_2'$ and a relabeling
$\sigma:V_1\rightarrow V_2'$ such that $V_1\cap V_2'=F_1$,
$\sigma(F_2)=F_1$, and $G_1|_{F_1}$ and $G_2|_{F_2}$ are complete
graphs on $t$ vertices.

Let $A$ be the set of pairs $(C,W)$ such that $W\subset V_1\cup V_2'$,
$|W|=k$ and $C$ is a connected component of $G|_W$.  Let
$$A_1 =\{(C,W)\in A: C\cap V_1\ne\emptyset\},\quad
A_2 =\{(C,W)\in A: C\cap V_2'\ne\emptyset\}.$$
Then $c_k(G)=|A|=|A_1|+|A_2|-|A_1\cap A_2|$.  It is sufficient to show
that $|A_1|=\sum_{i=0}^k c_i (G_1) \binom{n_2 - t}{k-i}$,
$|A_2|=\sum_{i=0}^k c_i (G_2) \binom{n_1 - t}{k-i}$ and $|A_1\cap
A_2|= \binom{n_1 + n_2 -t}{k} - \binom{n_1 + n_2 -2t}{k}$.  

Let $B_1$ be the set of triples $(C_1,W_1,X)$ such that $W_1\subset
V_1$, $X\subset V_2'\setminus V_1$, $|X|+|W_1|=k$ and $C_1$ is a
connected component of $G_1|_{W_1}$. Let $\phi_1:A_1\rightarrow B_1$
be the map defined by $\phi_1(C,W)=(C\cap V_1,W\cap V_1,W\setminus
V_1)$.  Then $\phi_1$ has the inverse map defined as follows. For a
triple $(C_1,W_1,X)\in B_1$, $\phi_1^{-1}(C_1,W_1,X)=(C,W)$, where
$W=W_1\cup X$ and $C$ is the connected component of $G|_{W}$
containing $C_1$. Thus $\phi_1$ is a bijection and we get
$|A_1|=|B_1|=\sum_{i=0}^k c_i (G_1) \binom{n_2 - t}{k-i}$.  Similarly
we get $|A_2|=\sum_{i=0}^k c_i (G_2) \binom{n_1 - t}{k-i}$.

Now let $B=\{W\subset V_1\cup V_2':W\cap F_1\ne \emptyset\}$.  Let
$\psi:A_1\cap A_2\rightarrow B$ be the map defined by $\psi(C,W)=W$.
We have the inverse map $\psi^{-1}$ as follows. For $W\in B$,
$\psi^{-1}(W)=(C,W)$, where $C$ is the connected component of $G|_W$
containing $W\cap F_1$, which is guaranteed to exist since
$G|_{F_1}=G_1|_{F_1}$ is a complete graph.  Thus $\psi$ is a
bijection, and we get $|A_1\cap A_2|= |B|=\binom{n_1 + n_2 -t}{k} -
\binom{n_1 + n_2 -2t}{k}$.
\end{proof}

\begin{thm}\label{thm:connectedsum}
Let $G_1$ and $G_2$ be graphs with $n_1$ and $n_2$ vertices
respectively. Let $t$ be a positive integer and let $G$ be a
$t$-connected sum of $G_1$ and $G_2$. Then
$$b_k (G) = \sum_{i=0}^k \left( b_i (G_1) \binom{n_2 - t}{k-i} + b_i
(G_2) \binom{n_1 - t}{k-i}\right) + \binom{n_1 + n_2 -2t}{k}.$$
\end{thm}
\begin{proof}
Since $c_k(G)=b_k(G)+\binom{n_1+n_2-t}{k}$,
$c_i(G_1)=b_i(G_1)+\binom{n_1}{i}$ and
$c_i(G_2)=b_i(G_2)+\binom{n_2}{i}$, by \autoref{lem:connectedsum}, it
is sufficient to show that
$$2\binom{n_1+n_2-t}{k} =
\sum_{i=0}^k \left( \binom{n_1}{i}\binom{n_2 - t}{k-i} + 
\binom{n_2}{i}\binom{n_1 - t}{k-i}\right),$$
which is immediate from the identity
$\sum_{i=0}^k\binom{a}{i}\binom{b}{k-i}=\binom{a+b}{k}$.
\end{proof}

Recall that a $t$-connected sum $G$ of two graphs depends on the
choice of vertices of each graph and the identification of the chosen
vertices. However, \autoref{thm:connectedsum} says that $b_k(G)$ does
not depend on them. Thus we get the following important property of a
$t$-connected sum of graphs.

\begin{cor}\label{thm:independent}
Let $t$ be a positive integer and let $G$ be a $t$-connected sum of
graphs $G_1,G_2,\dots,G_n$.  If $H$ is also a $t$-connected sum of
$G_1,G_2,\dots,G_n$, then $b_k(G)=b_k(H)$ for all $k$.
\end{cor}

Using \autoref{thm:tcon}, we get a formula for the special graded
Betti number of a $t$-connected sum of two simplicial complexes for
$t\geq3$.
\begin{cor}
Let $\Delta_1$ and $\Delta_2$ be simplicial complexes on $V_1$ and
$V_2$ respectively with $|V_1|=n_1$ and $|V_2|=n_2$. Let $t$ be a
positive integer and let $\Delta$ be a $t$-connected sum of $\Delta_1$
and $\Delta_2$. If $t\geq3$ then
$$ b_k (\Delta) = \sum_{i=0}^k \left( b_i (\Delta_1) \binom{n_2 -
  t}{k-i} + b_i (\Delta_2) \binom{n_1 - t}{k-i}\right) + \binom{n_1 +
  n_2 -2t}{k}.$$
\end{cor}

For an integer $n$, let $K_n$ denote a complete graph with $n$
vertices.

Let $G$ be a graph with vertex set $V$. If $H$ is a $t$-connected sum
of $G$ and $K_{t+1}$ then $H$ is a graph obtained from $G$ by adding a
new vertex $v$ connected to all vertices in $W$ for some $W\subset V$
such that $G|_W$ is isomorphic to $K_t$. Thus $H$ is determined by
choosing such a subset $W\subset V$. Using this observation, we get
the following lemma.

\begin{thm}\label{lem:complete}
Let $t$ be a positive integer.  Let $G$ be a $t$-connected sum of $n$
$K_{t+1}$'s. Then
$$b_k(G)=(k-1)\binom{n}{k}.$$
\end{thm}
\begin{proof}
We construct a sequence of graphs $H_1,\ldots,H_n$ as follows. Let
$H_1$ be the complete graph with vertex set
$\{v_1,v_2,\ldots,v_{t+1}\}$.  For $i\geq2$, let $H_i$ be the graph
obtained from $H_{i-1}$ by adding a new vertex $v_{t+i}$ connected to
all vertices in $\{v_1,v_2,\ldots,v_t\}$. Then $H_n$ is a
$t$-connected sum of $n$ $K_{t+1}$'s, and we have $b_k(G)=b_k(H_n)$ by
\autoref{thm:independent}. In $H_n$, the vertex $v_i$ is connected to
all the other vertices for $i\leq t$, and $v_j$ and $v_{j'}$ are not
connected to each other for all $t+1\leq j,j'\leq t+n$. Thus
$b_k(H_n)=(k-1)\binom{n}{k}$.
\end{proof}

Observe that every tree with $n+1$ vertices is a $1$-connected sum of
$n$ $K_2$'s.  Thus we get the following nontrivial property of trees
which was observed by Bruns and Hibi \cite{Bruns1998}.

\begin{cor}\cite[Example 2.1. (b)]{Bruns1998}\label{thm:tree}
Let $T$ be a tree with $n+1$ vertices. Then $b_k(T)$ does not depend
on the specific tree $T$.  We have
$$b_k(T)=(k-1)\binom{n}{k}.$$
\end{cor}

\begin{cor}\cite[Example 2.1. (c)]{Bruns1998}\label{thm:gon}
Let $G$ be an $n$-gon. If $k=n$ then $b_k(G)=0$; otherwise,
$$b_k(G) = \frac{n(k-1)}{n-k}\binom{n-2}{k}.$$
\end{cor}
\begin{proof}
It is clear for $k=n$. Assume $k<n$. Let
$V=\{v_1,\ldots,v_{n}\}$ be the vertex set of $G$.
Then 
\begin{align*}
(n-k)\cdot b_k(G) &= \sum_{\substack{W\subset V\\|W|=k}}
(\nc(G|_W)-1)\sum_{v\in V\setminus W} 1\\
&=\sum_{v\in V}\sum_{\substack{W\subset V\setminus\{v\}\\|W|=k}}
(\nc(G|_W)-1)\\
&=\sum_{v\in V} b_k(G|_{V\setminus\{v\}}).
\end{align*}
Since each $G|_{V\setminus\{v\}}$ is a tree with $n-1$ vertices, we
are done by \autoref{thm:tree}.
\end{proof}

\begin{remark}
Bruns and Hibi \cite{Bruns1998} obtained \autoref{thm:tree} and
\autoref{thm:gon} by showing that if $\Delta$ is a tree (or an
$n$-gon), considered as a $1$-dimensional simplicial complex, then
$\k[\Delta]$ has a pure resolution. Since $\k[\Delta]$ is
Cohen-Macaulay and it has a pure resolution, the Betti numbers are
determined by its type (c.f. \cite{Bruns1993}).
\end{remark}

Now we can prove \eqref{eq:hibi}.  Note that, for $d\geq3$, if $P$ is
a $d$-dimensional simplicial polytope and $Q$ is a simplicial polytope
obtained from $P$ by attaching a $d$-dimensional simplex $S$ to a
facet of $P$, then $\Delta(Q)$ is a $d$-connected sum of $\Delta(P)$
and $\Delta(S)$, and thus the $1$-skeleton $G(\Delta(Q))$ is a
$d$-connected sum of $G(\Delta(P))$ and $K_{d+1}$. Hence the
$1$-skeleton of the boundary complex of a $d$-dimensional stacked
polytope is a $d$-connected sum of $K_{d+1}$'s.

\begin{thm}
Let $P$ be a $d$-dimensional stacked polytope with $n$ vertices.  If
$d\geq3$, then $$b_k(\Delta(P)) = (k-1)\binom{n-d}{k}.$$ 
If $d=2$, then
$$b_k(\Delta(P)) = \left\{
\begin{array}{ll}
  0, & \mbox{if $k=n$,}\\
\frac{n(k-1)}{n-k}\binom{n-2}{k}, & \mbox{otherwise.}  
\end{array}\right.$$
\end{thm}
\begin{proof}
Assume $d\geq3$.  Then the $1$-skeleton $G(\Delta(P))$ is a
$d$-connected sum of $n-d$ $K_{d+1}$'s.  Thus by
\autoref{lem:complete}, we get
$b_k(\Delta(P))=b_k(G(\Delta(P)))=(k-1)\binom{n-d}{k}$.

Now assume $d=2$.  Then $G(\Delta(P))$ is an $n$-gon. Thus by
\autoref{thm:gon} we are done.
\end{proof}

\bibliographystyle{plain}

\end{document}